\documentclass{article}

\newcommand{\be}{\begin{equation}}
\newcommand{\ee}{\end{equation}}
\newcommand{\bea}{\begin{eqnarray}}
\newcommand{\eea}{\end{eqnarray}}
\newcommand{\barray}{\begin{array}}
\newcommand{\earray}{\end{array}}
\newcommand{\pa}{\partial}
\newcommand{\nn}{\nonumber}
\newcommand{\bitem}{\begin{itemize}}
\newcommand{\eitem}{\end{itemize}}
\newtheorem{teo}{Theorem}[section]
\newcommand{\bt}{\begin{teo}}
\newcommand{\et}{\end{teo}}
\newtheorem{Def}{Definition}[section]
\newcommand{\bd}{\begin{Def}}
\newcommand{\ed}{\end{Def}}
\newtheorem{lem}{Lemma}[section]
\newcommand{\bl}{\begin{lem}}
\newcommand{\el}{\end{lem}}
\newtheorem{prop}{Proposition}[section]
\newcommand{\bp}{\begin{prop}}
\newcommand{\ep}{\end{prop}}
\newtheorem{cor}{Corollary}[section]
\newcommand{\bc}{\begin{cor}}
\newcommand{\ec}{\end{cor}}
\newtheorem{ex}{Example}[section]
\newcommand{\bex}{\begin{ex}}
\newcommand{\eex}{\end{ex}}
\newtheorem{rem}{Remark}[section]
\newcommand{\br}{\begin{rem}}
\newcommand{\er}{\end{rem}}

\catcode `\@=11
\@addtoreset{equation}{section}

\begin{document}

\begin{center}
{\Large \textbf{On integrability of
the equations for nonsingular pairs of compatible flat
metrics\footnote{This work was
supported by the Alexander von Humboldt Foundation (Germany),
the Russian Foundation for Basic Research (project nos. 99--01--00010 and
 96--15--96027) and
INTAS (project no. 96--0770).}}}
\end{center}

\medskip

\begin{center}
{\large {O. I. Mokhov}}
\end{center}

\bigskip

\section{Introduction. Basic definitions} \label{vved}

In this paper, we deal with the problem of description of
nonsingular pairs of compatible flat metrics
 for the general $N$-component case.
We describe the scheme of
the integrating the nonlinear equations describing 
nonsingular pairs of compatible flat metrics (or, in other words,
nonsingular flat pencils of metrics).
This scheme was announced in our previous paper \cite{1}.
It is based on the reducing this problem to
a special
reduction of the Lam\'{e} equations and
the using the Zakharov method of differential
reductions \cite{za}
in the dressing method (a version of the inverse scattering
method).

We shall use both contravariant metrics $g^{ij} (u)$ with upper
indices,
where $u = (u^1,...,u^N)$ are local coordinates, $1 \leq i, j \leq N$,
and covariant metrics
$g_{ij}(u)$ with lower indices,
$g^{is} (u) g_{sj} (u) = \delta^i_j.$
The indices of coefficients of the Levi--Civita connections
 $\Gamma^i_{jk} (u)$
(the Riemannian connections generated by the corresponding
metrics) and tensors of Riemannian curvature $R^i_{jkl} (u)$
are raised and lowered by the metrics corresponding to them:
$$\Gamma^{ij}_k (u) = g^{is} (u) \Gamma^j_{sk} (u),
 \ \ \ \Gamma^i_{jk} (u) = {1 \over 2} g^{is} (u) \left (
{\pa g_{sk} \over \pa u^j} + {\pa g_{js} \over \pa u^k} -
{\pa g_{jk} \over \pa u^s} \right ),$$
$$R^{ij}_{kl} (u) = g^{is} (u) R^j_{skl} (u), \ \ \
R^i_{jkl} (u) = - {\pa \Gamma^i_{jl} \over \pa u^k}
+ {\pa \Gamma^i_{jk} \over \pa u^l} -
\Gamma^i_{pk} (u) \Gamma^p_{jl} (u)
+ \Gamma^i_{pl} (u) \Gamma^p_{jk} (u).$$

\bd [\cite{[31]}, \cite{[30]}]  \label{d}
Two contravariant flat metrics
$g_1^{ij} (u)$ and $g_2^{ij} (u)$ are called compatible
if any linear combination of these metrics
\be
g^{ij} (u) = \lambda_1 g_1^{ij} (u) + \lambda_2 g_2^{ij} (u),
\label{comb0}
\ee
where $\lambda_1$ and $\lambda_2$ are arbitrary constants such that
$\det ( g^{ij} (u) ) \not\equiv 0$,
is also a flat metric and
the coefficients of the corresponding Levi--Civita connections
are related by the same linear formula:
\be
\Gamma^{ij}_k (u) = \lambda_1 \Gamma^{ij}_{1, k} (u) +
\lambda_2 \Gamma^{ij}_{2, k} (u). \label{sv0}
\ee
We shall also say in this case that the metrics
$g_1^{ij} (u)$ and $g_2^{ij} (u)$ form a flat pencil.
\ed

\bd [\cite{1}] \label{dd}
Two contravariant metrics
$g_1^{ij} (u)$ and $g_2^{ij} (u)$ of constant Riemannian curvature
$K_1$ and $K_2$, respectively, are called compatible
if any linear combination of these metrics
\be
g^{ij} (u) = \lambda_1 g_1^{ij} (u) + \lambda_2 g_2^{ij} (u),
\label{comb00}
\ee
where $\lambda_1$ and $\lambda_2$ are arbitrary constants such that
$\det ( g^{ij} (u) ) \not\equiv 0$,
is a metric of constant Riemannian curvature
$\lambda_1 K_1 + \lambda_2 K_2$  and
the coefficients of the corresponding Levi--Civita connections
are related by the same linear formula:
\be
\Gamma^{ij}_k (u) = \lambda_1 \Gamma^{ij}_{1, k} (u) +
\lambda_2 \Gamma^{ij}_{2, k} (u). \label{sv00}
\ee
We shall also say in this case that the metrics
$g_1^{ij} (u)$ and $g_2^{ij} (u)$ form a pencil of
metrics of constant Riemannian curvature.
\ed

\bd [\cite{1}] \label{d1}
Two Riemannian or pseudo-Riemannian contravariant metrics
$g_1^{ij} (u)$ and $g_2^{ij} (u)$ are called compatible if
for any linear combination of these metrics
\be
g^{ij} (u) = \lambda_1 g_1^{ij} (u) + \lambda_2 g_2^{ij} (u),
\label{comb}
\ee
where $\lambda_1$ and $\lambda_2$ are arbitrary constants
such that $\det ( g^{ij} (u) ) \not\equiv 0$,
the coefficients of the corresponding Levi--Civita connections
and the components of the corresponding tensors of
Riemannian curvature are related by the same linear formula:
\be
\Gamma^{ij}_k (u) = \lambda_1 \Gamma^{ij}_{1, k} (u) +
\lambda_2 \Gamma^{ij}_{2, k} (u), \label{sv}
\ee
\be
R^{ij}_{kl} (u) = \lambda_1 R^{ij}_{1, kl} (u)
+ \lambda_2 R^{ij}_{2, kl} (u).  \label{kr}
\ee
We shall also say in this case that the metrics
$g_1^{ij} (u)$ and $g_2^{ij} (u)$ form a pencil of metrics.
\ed

\bd [\cite{1}] \label{d2}
Two Riemannian or pseudo-Riemannian contravariant metrics
$g_1^{ij} (u)$ and $g_2^{ij} (u)$ are called almost compatible
if for any linear combination of these metrics (\ref{comb})
relation (\ref{sv}) is fulfilled.
\ed

\bd
Two Riemannian or pseudo-Riemannian metrics
$g_1^{ij} (u)$ and $g_2^{ij} (u)$ are called nonsingular pair
of metrics if the eigenvalues of this pair of metrics, that is,
the roots of the equation
\be
\det ( g_1^{ij} (u) -  \lambda g_2^{ij} (u)) =0,
\ee
are distinct.
\ed

These definitions are motivated by the theory of
compatible Poisson brackets of hydrodynamic type.
In the case if the metrics $g_1^{ij} (u)$ and $g_2^{ij}(u)$ are
flat, that is,
 $R^i_{1, jkl} (u) = R^i_{2, jkl} (u) = 0,$
relation (\ref{kr}) is equivalent to the condition that
an arbitrary linear combination of the flat metrics
$g_1^{ij} (u)$ and $g_2^{ij}(u)$ is also a flat metric
and Definition \ref{d1} is equivalent to the well-known
definition of a flat pencil of metrics or, in other words,
a compatible pair of local nondegenerate
Poisson structures of hydrodynamic type
\cite{[31]} (see also \cite{[30]}--\cite{mokh2}).
In the case if the metrics $g_1^{ij} (u)$ and $g_2^{ij}(u)$ are
metrics of constant Riemannian curvature $K_1$ and $K_2$,
respectively, that is,
$$R^{ij}_{1, kl} (u) = K_1 (\delta^i_k \delta^j_l -
\delta^i_l  \delta^j_k), \ \ \
R^{ij}_{2, kl} (u) = K_2 (\delta^i_k \delta^j_l -
\delta^i_l  \delta^j_k),$$
relation (\ref{kr}) gives the condition that an arbitrary
linear combination of the metrics
 $g_1^{ij} (u)$ and $g_2^{ij}(u)$  (\ref{comb})
is a metric of
constant Riemannian curvature
$\lambda_1 K_1 +
\lambda_2 K_2$ and Definition \ref{d1} is
equivalent to Definition \ref{d} of a pencil of
metrics of constant Riemannian curvature or, in other words,
a compatible pair of the corresponding nonlocal Poisson structures
of hydrodynamic type which were introduced and studied by the author
and Ferapontov in \cite{[38]}.
Compatible metrics of more general type correspond
to compatible pairs of nonlocal Poisson structures
of hydrodynamic type which were introduced and studied by
Ferapontov in \cite{[40]}. They arise, for example,
if we shall use a recursion operator generated by a pair
of compatible Poisson structures of hydrodynamic type and
determining, as is well-known, an infinite sequence of
corresponding Poisson structures.

\section{Compatible local Poisson structures of
\\ hydrodynamic type (a brief survey)}

The local homogeneous Poisson bracket of the first order, that is,
the Poisson bracket of the form
\be
\{ u^i (x), u^j (y) \} =
g^{ij} (u(x))\, \delta_x (x-y) + b^{ij}_k (u(x)) \, u^k_x \,
\delta (x-y), \label{(1.1)}
\ee
where $u^1,...,u^N$ are local coordinates on a certain smooth
$N$-dimensional manifold $M$, is called a
{\it local Poisson structure of hydrodynamic type} or
{\it Dubrovin--Novikov structure} \cite{[1]}.
Here, $u^i(x),\ 1 \leq i \leq N,$ are functions (fields) of
a single independent variable $x$, and
the coefficients $g^{ij}(u)$ and $b^{ij}_k (u)$ of bracket (\ref{(1.1)})
are smooth functions on $M$.

In other words,
for arbitrary functionals $I[u]$ and $J[u]$
on the space of fields $u^i(x), \ 1 \leq i \leq N,$ a bracket of the form
\be
\{ I,J \} = \int {\delta I \over \delta u^i(x)}
\biggl ( g^{ij}(u(x)) {d \over dx} + b^{ij}_k (u(x))\, u^k_x \biggr )
{\delta J \over \delta u^j(x)} dx
\label{(1.2)}
\ee
is defined and it is required that this bracket is a Poisson bracket,
that is, it is skew-symmetric:
\be
\{ I, J \} = - \{ J, I \}, \label{skew}
\ee
and satisfies the Jocobi identity
\be
\{ \{ I, J \}, K \} + \{ \{ J, K \}, I \} + \{ \{ K, I \}, J \} =0
\label{jacobi}
\ee
for arbitrary functionals $I[u]$, $J[u]$ and $K[u]$.
The skew-symmetry (\ref{skew}) and the Jacobi identity
(\ref{jacobi}) impose very strict conditions on the coefficients
$g^{ij}(u)$ and $b^{ij}_k (u)$ of bracket (\ref{(1.2)}) (these
conditions will be considered below).

The local Poisson structures of hydrodynamic type (\ref{(1.1)})
were introduced and studied by Dubrovin and Novikov
in \cite{[1]}. In this paper, they proposed
a general Hamiltonian approach to the so-called
{\it homogeneous systems of hydrodynamic type}, that is, to
evolutionary quasilinear systems of first-order partial differential
equations
\be
u^i_t = V^i_j (u)\, u^j_x
\label{(1.5)}
\ee
that corresponds to structures (\ref{(1.1)}).

This Hamiltonian approach was motivated by the study
of the equations of Euler hydrodynamics and the
Whitham averaging equations that describe the evolution of
slowly modulated multiphase solutions of partial
differential equations \cite{[3]}.

Local bracket (\ref{(1.2)}) is called
{\it nondegenerate} if
$\det (g^{ij} (u)) \not\equiv 0$.
For general nondegenerate brackets of form (\ref{(1.2)}),
Dubrovin and Novikov proved the following important theorem.

\bt [Dubrovin, Novikov \cite{[1]}] \label{dn}
If $\det (g^{ij} (u)) \not\equiv 0$, then bracket (\ref{(1.2)})
is a Poisson bracket, that is, it is skew-symmetric and satisfies
the Jacobi identity if and only if
\bitem
\item [(1)] $g^{ij} (u)$ is an arbitrary flat
pseudo-Riemannian contravariant metric (a metric of
zero Riemannian curvature),

\item [(2)] $b^{ij}_k (u) = - g^{is} (u) \Gamma ^j_{sk} (u),$ where
$\Gamma^j_{sk} (u)$ is the Riemannian connection generated by
the contravariant metric $g^{ij} (u)$
(the Levi--Civita connection).
\eitem
\et

Consequently, for any local nondegenerate Poisson structure of
hydrodynamic type, there always exist local coordinates
$v^1,...,v^N$ (flat coordinates of the metric $g^{ij}(u)$) in which
the coefficients of the brackets are constant:
\be
\widetilde g^{ij} (v)
= \eta^{ij} = {\rm \ const}, \ \
\widetilde \Gamma^i_{jk} (v) = 0, \ \
\widetilde b^{ij}_k (v) =0,
\ee
that is, the bracket has the constant form
\be
\{ I,J \} = \int {\delta I \over \delta v^i(x)}
 \eta^{ij} {d \over dx}
{\delta J \over \delta v^j(x)} dx,
\label{(1.6)}
\ee
where $(\eta^{ij})$ is a nondegenerate symmetric constant matrix:
\be
\eta^{ij} = \eta^{ji}, \ \ \eta^{ij} = {\rm const},
\  \ \det \, (\eta^{ij}) \neq 0.\nn
\ee

On the other hand, as early as 1978, Magri proposed
a bi-Hamiltonian approach to the integration of nonlinear
systems \cite{[11]}. This approach demonstrated
that the integrability is closely related to the
bi-Hamiltonian property, that is, to the property
of a system to have two compatible Hamiltonian
representations. As was shown by Magri in \cite{[11]},
compatible Poisson brackets generate integrable hierarchies of
systems of differential equations.
Therefore, the description of compatible
Poisson structures is very urgent
and important problem in the theory of
integrable systems.
In particular, for a system, the bi-Hamiltonian property
generates recurrent relations for the conservation laws of
this system.

Beginning from \cite{[11]}, quite extensive literature
(see, for example,
\cite{[12]}--\cite{[16]}
and the necessary references therein) has been
devoted to the bi-Hamiltonian approach and to
the construction of compatible Poisson structures
for many specific important equations of
mathematical physics and field theory.
As far as the problem of description
of sufficiently wide classes
of compatible Poisson structures of defined special types is
concerned, apparently the first such statement was considered
in \cite{[17]}, \cite{[18]} (see also \cite{[19]}, \cite{[20]}).
In those papers, the present author posed and completely solved
the problem of description of all compatible local scalar first-order
and third-order Poisson brackets, that is, all Poisson brackets given
by arbitrary scalar first-order and third-order ordinary
differential
operators.
These brackets generalize the well-known compatible pair of
the Gardner--Zakharov--Faddeev bracket
\cite{[21]}, \cite{[22]}
(first-order bracket) and the Magri bracket \cite{[11]}
(third-order bracket) for the Korteweg--de Vries equation.

In the case of homogeneous systems of hydrodynamic type, many integrable
systems possess compatible Poisson structures of hydrodynamic type.
The problems of description of these structures for particular systems and
numerous examples were considered in many papers
(see, for example, \cite{[23]}--\cite{[29]}).
In particular, in \cite{[23]}
Nutku studied a special class of compatible
two-component Poisson structures of hydrodynamic type and
the related bi-Hamiltonian hydrodynamic systems.
In \cite{[28]} Ferapontov classified all two-component
homogeneous systems of hydrodynamic type possessing
three compatible local Poisson structures of
hydrodynamic type.

In the general form, the problem of description of
flat pencil of metrics (or, in other words,
compatible nondegenerate local Poisson structures of hydrodynamic type)
was considered by Dubrovin in
\cite{[31]}, \cite{[30]}
in connection with the construction of important examples of
such flat pencils of metrics, generated by natural pairs of
flat metrics on the spaces of orbits of Coxeter groups and on other
Frobenius manifolds and associated with the
corresponding quasi-homogeneous solutions of
the associativity equations.
In the theory of Frobenius manifolds introduced and studied by
Dubrovin
\cite{[31]}, \cite{[30]}
(they correspond to two-dimensional topological field theories),
a key role is played by flat pencils of metrics, possessing
a number of special additional (and very strict) properties
(they satisfy the so-called quasi-homogeneity property).
In addition, in \cite{[33]} Dubrovin proved that the theory of
Frobenius manifolds is equivalent to the theory
quasi-homogeneous compatible nondegenerate Poisson structures of
hydrodynamic type. The general problem of compatible
nondegenerate local Poisson structures was also considered by Ferapontov
in \cite{[34]}.

The author's papers
\cite{[35]}--\cite{mokh2}
are devoted to the general problem of classification of
local Poisson structures of hydrodynamic type, to integrable
nonlinear systems which describe such
compatible Poisson structures and to special reductions connected
with the associativity equations.

\bd [Magri \cite{[11]}] \label{dm}
Two Poisson brackets $\{ \ , \ \}_1$
and $\{ \ , \ \}_2$ are called {\it compatible} if
an arbitrary linear combination of these Poisson brackets
\be
\{ \ , \ \} = \lambda_1 \, \{ \ , \ \}_1 +
\lambda_2 \, \{ \ , \ \}_2, \label{magri}
\ee
where $\lambda_1$ and $\lambda_2$ are arbitrary constants,
is also always a Poisson bracket.
In this case, one can say also that the brackets $\{ \ , \ \}_1$ and
$\{ \ , \ \}_2$ form a pencil of Poisson brackets.
\ed

Correspondingly, the problem of description
of compatible nondegenerate local Poisson structures
of hydrodynamic type is pure differential-geometric
problem of description of flat pencils of metrics
(see \cite{[31]}, \cite{[30]}).

In \cite{[31]}, \cite{[30]}
Dubrovin presented all the tensor relations
for the general flat pencils of metrics.
First, we introduce the necessary notation.
Let $\nabla_1$ and $\nabla_2$
be the operators of covariant differentiation given by the Levi--Civita
connections $\Gamma^{ij}_{1,k} (u)$
and $\Gamma^{ij}_{2,k} (u)$, generated by the metrics $g^{ij}_1 (u)$
and $g^{ij}_2 (u)$, respectively. The indices
of the covariant differentials are raised and lowered by
the
corresponding metrics: $\nabla^i_1= g^{is}_1 (u) \nabla_{1,s}$,
$\nabla^i_2=g^{is}_2 (u) \nabla_{2,s}$.
Consider the tensor
\be
\Delta ^{ijk} (u) = g^{is}_1 (u) g^{jp}_2 (u)
\left (\Gamma^k_{2, ps} (u)
- \Gamma^k_{1, ps} (u) \right ),  \label{(2.3)}
\ee
introduced by Dubrovin in \cite{[31]}, \cite{[30]}.

\bt [Dubrovin \cite{[31]}, \cite{[30]}] \label{dub1}
 If metrics
$g^{ij}_1 (u)$ and $g^{ij}_2 (u)$ form a flat pencil,
then there exists a vector field $f^i (u)$ such that the tensor
$\Delta ^{ijk} (u)$ and the metric $g^{ij}_1 (u)$ have the form
\be
\Delta ^{ijk} (u) = \nabla^i_2 \nabla^j_2 f^k (u),
\label{(2.4)}
\ee
\be
g^{ij}_1 (u) = \nabla^i_2 f^j (u) + \nabla^j_2 f^i (u) +
c g^{ij}_2 (u), \label{(2.5)}
\ee
where $c$ is a certain constant,
and the vector field $f^i (u)$ satisfies the equations
\be
\Delta^{ij}_s (u) \Delta^{sk}_l (u) =
\Delta^{ik}_s (u) \Delta^{sj}_l (u), \label{(2.6)}
\ee
where
\be
\Delta^{ij}_k (u) =g_{2,ks} (u) \Delta ^{sij} (u)
= \nabla_{2,k} \nabla^i_2 f^j (u), \label{(2.7)}
\ee
and
\be
(g^{is}_1 (u) g^{jp}_2 (u) - g^{is}_2 (u) g^{jp}_1 (u))
\nabla_{2,s} \nabla_{2,p} f^k (u) =0. \label{(2.8)}
\ee
Conversely, for the flat metric $g^{ij}_2 (u)$ and the vector field
$f^i (u)$ that is a solution of the system of equations
(\ref{(2.6)}) and (\ref{(2.8)}),
the metrics $g^{ij}_2 (u)$ and (\ref{(2.5)}) form a flat pencil.
\et

The proof of this theorem immediately follows from the relations that
are equivalent to the fact that the metrics $g^{ij}_1 (u)$
and $g^{ij}_2 (u)$ form a flat pencil and are considered in
flat coordinates of the metric $g^{ij}_2 (u)$ \cite{[31]}, \cite{[30]}.

In my paper \cite{[35]}, an explicit and simple {\em criterion}
of compatibility for two Poisson structures of
hydrodynamic type is formulated, that is, it is shown
what explicit form is sufficient and necessary for the
Poisson structures of hydrodynamic type to be
compatible.

For the moment, we are able to formulate such explicit general criterion
only namely in terms of Poisson structures but not in terms of
metrics as in Theorem \ref{dub1}.

\bl [\cite{[35]}]  \label{lem1}
{\bf (An explicit criterion of
compatibility for Poisson structures of
hydrodynamic type)}
Any local Poisson structure of hydrodynamic type
$\{ I, J \}_2$  is compatible with the constant nondegenerate
Poisson bracket (\ref{(1.6)}) if and only if it has the form
\be
\{ I, J \}_2 =
\int {\delta I \over \delta v^i (x)} \biggl (
\biggl ( \eta^{is} {\partial h^j \over \partial v^s}
+ \eta^{js} {\partial h^i \over \partial v^s} \biggr )
{d \over dx} + \eta^{is}
{\partial^2 h^j \over \partial v^s \partial v^k}
v^k_x \biggr ) {\delta J \over \delta v^j (x)} dx, \label{(2.9)}
\ee
where $h^i (v), \ 1 \leq i \leq N,$ are smooth
functions defined on a certain neighbourhood.
\el

We do not require in Lemma \ref{lem1} that
the Poisson structure of hydrodynamic type
$\{ I, J \}_2$ is nondegenerate.
Besides, it is important to note that this statement is local.

In 1995, in the paper \cite{[34]}, Ferapontov
proposed an approach to the problem on flat pencils of metrics,
which is motivated by the theory of recursion operators,
and formulated the following theorem as a criterion of
compatibility of nondegenerate local Poisson structures
of hydrodynamic type:

\bt [\cite{[34]}] \label{tfer}
Two local nondegenerate Poisson structures of hydrodynamic type
given by flat metrics $g_1^{ij}(u)$ and $g_2^{ij}(u)$
are compatible if and only if the Nijenhuis tensor of the affinor
$v^i_j (u) = g_1^{is} (u) g_{2, sj} (u)$ vanishes,
that is,
\be
N^k_{ij} (u) = v^s_i (u) {\pa v^k_j \over \pa u^s}
- v^s_j (u) {\pa v^k_i \over \pa u^s}  +
v^k_s (u) {\pa v^s_i \over \pa u^j} -
v^k_s {\pa v^s_j \over \pa u^i} =0.
\ee
\et

Besides, it is noted in the remark in \cite{[34]} that
if the spectrum of $v^i_j (u)$ is simple,
the vanishing of the Nijenhuis tensor implies
the existence of coordinates $R^1,...,R^N$ for which
all the objects $v^i_j (u)$, $g_1^{ij} (u)$, $g_2^{ij} (u)$
become diagonal.
Moreover, in these coordinates the
$i$th eigenvalue of $v^i_j (u)$ depends only on the
coordinate $R^i$.
In the case when all the eigenvalues are nonconstant,
they can be introduced as new coordinates.
In these new coordinates
$\tilde v^i_j (R) = {\rm diag \ } (R^1,...,R^N)$,
$\tilde g_2^{ij} (R) = {\rm diag \ } (g^1 (R),...,g^N (R))$,
$\tilde g_1^{ij} (R) = {\rm diag \ } (R^1 g^1 (R),..., R^N g^N (R))$.

Unfortunately, as is shown in \cite{1},
in the general case the Theorem \ref{tfer} is not true
and, correspondingly, it is not a criterion of
compatibility of flat metrics.
Generally speaking, compatibility of flat metrics does not follow
from the vanishing of the corresponding Nijenhuis tensor.
The corresponding counterexamples
were presented in \cite{1}.
In the general case, as it was shown in \cite{1},
the Theorem \ref{tfer} is actually a criterion of
almost compatibility of flat metrics that
does not guarantee compatibility of the corresponding
nondegenerate local Poisson structures of hydrodynamic type.
But if the spectrum of $v^i_j (u)$ is simple, that is, all the
eigenvalues are distinct, then
the Theorem \ref{tfer} is not only true but also can be
essentially  generalized for the case of arbitrary
compatible Riemannian or pseudo-Riemannian metrics,
in particular, for the especially important cases in
the theory of systems of hydrodynamiic type, namely,
the cases of metrics of constant Riemannian curvature
or the metrics generating the general nonlocal Poisson
structures of hydrodynamic type (see \cite{1}).

In particular, the following general theorem is proved in \cite{1}:
\bt  [\cite{1}]   \label{teom}
An arbitrary nonsingular pair of metrics is compatible
if and only if there exist local coordinates $u = (u^1,...,u^N)$
such that
$g^{ij}_2 (u) = g^i (u) \delta^{ij}$ and
$g^{ij}_1 (u) = f^i (u^i) g^i (u) \delta^{ij},$
where $f^i (u^i),$ $i=1,...,N,$ are arbitrary functions
of single variables (of course, in the case of nonsingular pair of metrics,
these functions are not equal to each other if they are constants and they
are not equal identically to zero).
\et

In this paper, we consider only the
case of nonsingular pairs of flat metrics.
In this case the approach of Ferapontov and Theorem \ref{tfer}
are absolutely correct.

\section{Equations for nonsingular pairs of \\
compatible flat metrics}

Let us consider here the problem on nonsingular pairs
of compatible flat metrics. It follows from Theorem \ref{tfer} and
Theorem \ref{teom}
that it is
sufficient to classify flat metrics of the form
$g_2^{ij} (u) = g^i (u) \delta^{ij}$ and
$g_1^{ij} (u) = f^i (u^i) g^i (u) \delta^{ij},$
where $f^i (u^i),$ $i= 1,...,N,$ are arbitrary functions of
single variables.

The  problem of description of diagonal flat metrics,
that is, flat metrics
$g_2^{ij} (u) = g^i (u) \delta^{ij},$
is a classical problem of differential geometry.
This problem is equivalent to the problem
of description of curvilinear orthogonal coordinate systems in
a pseudo-Euclidean space and it was studied in detail and mainly solved
in the beginning of the 20th century (see \cite{da}).
Locally, such coordinate systems are determined by
$n(n-1)/2$ arbitrary functions of two variables (see \cite{ca},
\cite{bi}). Recently, Zakharov showed that the Lam\'{e}
equations describing curvilinear orthogonal coordinate systems
can be integrated by the inverse scattering method \cite{za}
(see also an algebraic-geometric approach in \cite{kri}).

\bt \label{syst}
Nonsingular pairs of compatible flat metrics are described by
the following integrable nonlinear
systems which are the special reductions of the Lam\'{e} equations:
\be
{\pa \beta_{ij} \over \pa u^k}
=\beta_{ik} \beta_{kj},\ \ \ i\neq j,\ \ i\neq k,\ \ j\neq k, \label{lam1}
\ee
\be
{\pa \beta_{ij} \over \pa u^i}+
{\pa \beta_{ji} \over \pa u^j}+\sum_{s\neq i,\
s\neq j} \beta_{si} \beta_{sj} =0,\ \ \ i\neq j, \label{lam2}
\ee
\be
\sqrt {f^i (u^i)} {\pa \left (\sqrt {f^i (u^i)} \beta_{ij}\right )
 \over \pa u^i}+
\sqrt {f^j (u^j)}
{\pa \left (\sqrt {f^j (u^j)}
\beta_{ji}\right ) \over \pa u^j}+\sum_{s\neq i,\
s\neq j} f^s (u^s) \beta_{si} \beta_{sj} =0,\ \ \ i\neq j, \label{lam3}
\ee
where $f^i (u^i),$ $i=1,...,N,$ are the given arbitrary functions
of single variables.
\et

The equations (\ref{lam1}) and (\ref{lam2}) are
the famous Lam\'{e} equations and the equation (\ref{lam3})
defines a nontrivial nonlinear differential reduction of
the Lam\'{e} equations. Such types of differential reductions
for Lam\'{e} equations
were also studied by Zakharov \cite{za} and the Zakharov method
can be applied successfully to our problem.

Consider the conditions of flatness for the diagonal metrics
$g^{ij}_2 (u) = g^i (u) \delta^{ij}$ and
$g^{ij}_1 (u) = f^i (u^i) g^i (u) \delta^{ij},$
where $f^i (u^i),$ $i=1,...,N,$ are arbitrary functions of
the given single variables (but these functions
are not equal to zero identically).

Recall that for any diagonal metric
$\Gamma^i_{jk} (u) =0$ if all the indices $i, j, k$
are distinct.
Correspondingly,
$R^{ij}_{kl} (u) = 0$
if all the indices $i, j, k, l $  are distinct.
Besides, as a result of the well-known symmetries of
the tensor of Riemannian curvature we have:
$$R^{ii}_{kl} (u) = R^{ij}_{kk} (u) =0,$$
$$R^{ij}_{il} (u) = -R^{ij}_{li} (u) = R^{ji}_{li} (u) =
- R^{ji}_{il} (u).$$

Thus, it is sufficient to consider
the condition $R^{ij}_{kl} (u) =0$ (the condition of the flatness for a
metric) only
for the following components of the tensor of Riemannian curvature:
 $R^{ij}_{il} (u)$,
where $i \neq j,$ $\ i \neq  l$.

For any diagonal metric
$g^{ij}_2 (u) = g^i (u) \delta^{ij}$   we have
$$\Gamma^i_{2, ik} (u) = \Gamma^i_{2, ki} (u) =
- {1 \over 2 g^i (u)} {\pa g^i \over \pa u^k}, \ \ \
{\rm \ for\  any\ } i, k; $$
$$\Gamma^i_{2, jj} (u) = {1 \over 2} {g^i (u) \over
(g^j (u))^2 } {\pa g^j \over \pa u^i},\ \ \ i \neq j.$$

\bea
&&
R^{ij}_{2, il} (u) = g^i (u) R^j_{2, iil} (u) = \nn\\
&&
g^i (u) \left ( - {\pa \Gamma^j_{2, il} \over \pa u^i}  +
{\pa \Gamma^j_{2, ii}  \over \pa u^l} - \sum_{s=1}^N
\Gamma^j_{2, si} (u) \Gamma^s_{2, il} (u) +
\sum_{s=1}^N \Gamma^j_{2, sl} (u) \Gamma^s_{2, ii} (u) \right ).
\eea

It is necessary to consider separately two different cases.

1) $j \neq l$.

\bea
&&
R^{ij}_{2, il} (u) =
g^i (u) \left (
{\pa \Gamma^j_{2, ii}  \over \pa u^l} -
\Gamma^j_{2, ii} (u) \Gamma^i_{2, il} (u) +
\Gamma^j_{2, jl} (u) \Gamma^j_{2, ii} (u)
 +
\Gamma^j_{2, ll} (u) \Gamma^l_{2, ii} (u) \right )
=  \nn\\
&&
{1 \over 2} g^i (u)  {\pa \over \pa u^l}
\left ( {g^j (u) \over (g^i (u))^2 } {\pa g^i \over \pa u^j }\right )
 + {1 \over 4} {g^j (u) \over (g^i (u))^2 }
{\pa g^i \over \pa u^j} {\pa g^i \over \pa u^l} \nn\\
&&
- {1 \over 4 g^i (u)}
{\pa g^i \over \pa u^j} {\pa g^j \over \pa u^l}
+ {1 \over 4} {g^j (u) \over g^i (u) g^l (u) }
{\pa g^l \over \pa u^j} {\pa g^i \over \pa u^l}=0 .\label{r1}
\eea

Introducing the standard classic notation
\bea
&&
g^i (u) = {1 \over (H_i (u))^2 },\ \ \ d\, s^2 = \sum_{i=1}^N
(H_i (u))^2 (d u^i)^2, \\
&&
\beta_{ik} (u) = {1 \over H_i (u)} {\pa H_k \over \pa u^i},\ \ \
i \neq k,
\eea
where $H_i (u)$ are the {\it Lam\'{e} coefficients} and
$\beta_{ik} (u)$ are the {rotation coefficients},
we obtain that equations (\ref{r1}) are equivalent to
the equations
\be
{\pa^2 H_i \over \pa u^j \pa u^k} =
{1 \over H_j (u)} {\pa H_i \over \pa u^j} {\pa H_j \over \pa u^k}
+ {1 \over H_k (u)} {\pa H_k \over \pa u^j} {\pa H_i \over \pa u^k}
\ee
or, equivalently, to equations (\ref{lam1}).

2) $j=l$.

\bea
&&
R^{ij}_{2, ij} (u) =
g^i (u) \left (
- {\pa \Gamma^j_{2, ij}  \over \pa u^i} +
{\pa \Gamma^j_{2, ii}  \over \pa u^j} -
\Gamma^j_{2, ii} (u) \Gamma^i_{2, ij} (u) - \right. \nn\\
&&
\left. \Gamma^j_{2, ji} (u) \Gamma^j_{2, ij} (u)
 +  \sum_{s=1}^N
\Gamma^j_{2, sj} (u) \Gamma^s_{2, ii} (u) \right )
=  \nn\\
&&
{1 \over 2} g^i (u)  {\pa \over \pa u^i}
\left ( {1 \over g^j (u) } {\pa g^j \over \pa u^i }\right )
+ {1 \over 2} g^i (u)  {\pa \over \pa u^j}
\left ( {g^j (u) \over (g^i (u))^2 } {\pa g^i \over \pa u^j }\right )
+ {1 \over 4} {g^j (u) \over (g^i (u))^2 }
{\pa g^i \over \pa u^j} {\pa g^i \over \pa u^j} \nn\\
&&
- {1 \over 4} {g^i (u) \over (g^j(u))^2}
{\pa g^j \over \pa u^i} {\pa g^j \over \pa u^i}
+ {1 \over 4 g^j (u)}
{\pa g^j \over \pa u^i} {\pa g^i \over \pa u^i} -
\sum_{s \neq  i} {1 \over 4} {g^s (u) \over g^i (u) g^j (u) }
{\pa g^j \over \pa u^s} {\pa g^i \over \pa u^s}=0.\label{r2}
\eea

Equations (\ref{r2}) are equivalent to the equations
\be
{\pa \over \pa u^i} \left ( {1 \over H_i (u)}
{\pa H_j  \over \pa u^i} \right ) +
{\pa \over \pa u^j} \left ( {1 \over H_j (u)}
{\pa H_i  \over \pa u^j} \right ) +
\sum_{s \neq i,\ s \neq j}
{1 \over (H_s (u))^2} {\pa H_i \over \pa u^s} {\pa H_j \over \pa u^s},
\ \ \ i \neq j,
\ee
or, equivalently, to equations (\ref{lam2}).

The condition that the metric
$g_1^{ij} (u) = f^i (u^i) g^i (u) \delta^{ij}$ is also flat
gives exactly $n(n-1)/2$ additional equations (\ref{lam3})
which are linear with respect to
the given functions $f^i (u^i)$.  Note that, in this case,
components (\ref{r1})
of tensor of Riemannian curvature
automatically vanish.
And the vanishing of components (\ref{r2}) gives
the corresponding $n(n-1)/2$ equations.

Actually, for the metric
$g^{ij}_1 (u) = f^i (u^i) g^i (u) \delta^{ij},$
we have
\be
\widetilde H_i (u) = {H_i (u)  \over \sqrt {f^i (u^i)}},\ \ \
\widetilde \beta_{ik} (u) = {\sqrt {f^i (u^i)} \over
\sqrt {f^k (u^k)}} \left ( {1 \over H_i (u)} {\pa H_k \over
\pa u^i}  \right ) =  {\sqrt {f^i (u^i)} \over
\sqrt {f^k (u^k)}}  \beta_{ik} (u).
\ee
Respectively,
equations (\ref{lam1}) are also fulfilled for
the rotation coefficients $\widetilde \beta_{ik} (u)$
and equations (\ref{lam2}) for them give equation
(\ref{lam3}).

 In particular, in the case
$N=2$ this completely solves the problem of description
of nonsingular pairs of compatible flat metrics \cite{1}.
In the next section we give their complete description.

\section{Two-component compatible flat metrics} \label{sect5}

We present here the complete description of
nonsingular pairs of two-component compatible
flat metrics \cite{1} (see also \cite{[35]}, \cite{mokh1},
\cite{mokh2}, where an integrable
four-component
 homogeneous system of hydrodynamic type,
describing all the two-component compatible flat metrics, was derived and
investigated).

It is shown above that for any nonsingular pair
of two-component compatible metrics $g^{ij}_1 (u)$
and $g^{ij}_2 (u)$ there always exist local coordinates
$u^1,...,u^N$ such that
\be
(g^{ij}_2 (u) ) = \left ( \begin{array}{cc}
{\varepsilon^1 \over (b^1 (u))^2 } & 0\\
0 & {\varepsilon^2 \over (b^2 (u))^2 }
\end{array} \right ), \ \ \ \ \
(g^{ij}_1 (u) ) = \left ( \begin{array}{cc}
{\varepsilon^1  f^1 (u^1) \over (b^1 (u))^2 } & 0\\
0 & {\varepsilon^2  f^2 (u^2) \over (b^2 (u))^2 }
\end{array} \right ), \label{metr}
\ee
where $\varepsilon^i = \pm 1, \ i=1, 2;$
$b^i (u)$ and $f^i (u^i), \ i=1, 2,$ are arbitrary nonzero functions
of the corresponding single variables.

\bl \label{moflat}
An arbitrary diagonal metric $g^{ij}_2 (u)$ (\ref{metr})
is flat if and only if the functions $b^i (u),$  $ i=1, 2,$
are solutions of the following linear system:
\be
{\pa b^2 \over \pa u^1} = \varepsilon^1
{\pa F \over \pa u^2} b^1 (u),\ \ \ \
{\pa b^1 \over \pa u^2} =
- \varepsilon^2 {\pa F \over \pa u^1} b^2 (u),
\label{sys}
\ee
where $F(u)$ is an arbitrary function.
\el

\bt [\cite{1}]
The metrics $g^{ij}_1 (u)$ and $g^{ij}_2 (u)$ (\ref{metr})
form a flat pencil of metrics if and only if
the functions $b^i (u), \ i=1, 2,$ are solutions of
the linear system (\ref{sys}), where the function $F(u)$ is
a solution of the following linear equation:
\be
2 {\pa^2 F \over \pa u^1 \pa u^2} (f^1 (u^1) -  f^2 (u^2))
+ {\pa F \over \pa u^2} {d f^1 (u^1)  \over d u^1} -
{\pa F \over \pa u^1} {d f^2 (u^2) \over d u^2} =0. \label{lequa}
\ee
\et

In the case, if the eigenvalues of the pair of the metrics
$g^{ij}_1 (u)$ and $g^{ij}_2 (u)$ are not only distinct but
also are not constants, we can always choose local coordinates
such that $f^1 (u^1) = u^1,$  $f^2 (u^2) = u^2$ (see also remark in
\cite{[34]}).
In this case, equation (\ref{lequa}) has the form
\be
2 {\pa^2 F \over \pa u^1 \pa u^2} (u^1 -  u^2)
+ {\pa F \over \pa u^2} -
{\pa F \over \pa u^1} = 0. \label{lequa2}
\ee

Let us continue this recurrent procedure for the metrics
$G^{ij}_{n+1} (u) = v^i_j (u) G^{ij}_n(u)$ with the help of
the affinor $v^i_j (u) = u^i \delta^i_j.$

\bt [\cite{1}]
Three metrics
\be
(G^{ij}_n (u) ) = \left ( \begin{array}{cc}
{\varepsilon^1  (u^1)^n \over (b^1 (u))^2 } & 0\\
0 & {\varepsilon^2   (u^2)^n \over (b^2 (u))^2 }
\end{array} \right ), \ \ \ \ \ n=0, 1, 2,\label{metr2}
\ee
form a flat pencil of metrics (pairwise compatible)
if and only if the functions
$b^i (u), \ i=1, 2,$ are solutions of the linear system (\ref{sys}),
where
\be
F(u) = c \ln (u^1 - u^2),
\ee
$c$ is an arbitrary constant.
Already the metric $G^{ij}_3 (u)$ is flat only in the
most trivial case, when $c =0,$
and, respectively, $b^1 = b^1 (u^1),$ $b^2 = b^2 (u^2)$.

The metric $G^{ij}_3 (u)$ is a metric of
nonzero constant Riemannian curvature $K \neq 0$
(in this case, the metrics $G^{ij}_n,$ $n=0, 1, 2, 3,$
form a pencil of metrics of constant Riemannian curvature)
if and only if
\be
(b^1 (u))^2 = (b^2 (u))^2 = {\varepsilon^2 \over 4 K}
(u^1 - u^2),\ \ \ \ \varepsilon^1 = - \varepsilon^2,
\ \ \ \ c = \pm {1 \over 2}.
\ee
\et

\section{Compatible flat metrics and the Zakharov
method of differential reductions}

Recall very briefly the Zakharov method of the integrating
the Lam\'{e} equations (\ref{lam1}) and (\ref{lam2}) \cite{za}.

We should choose a matrix function $F_{ij} (s, s', u)$ and
solve the integral equation
\be
K_{ij}(s, s', u)  = F_{ij}(s, s', u) + \int_s^{\infty}
\sum_l K_{il} (s, q, u) F_{lj} (q, s', u) dq.  \label{int}
\ee
Then we obtain a one-parameter family of solutions of
the Lam\'{e} equations by the
formula
\be
\beta_{ij} (s, u) = K_{ji} (s, s, u).  \label{resh}
\ee

In particular,
if $F_{ij} (s, s', u) = f_{ij} (s-u^i, s'-u^j),$
where $f_{ij} (x, y)$ is an arbitrary matrix function of two variables,
then formula (\ref{resh}) produces solutions of equations (\ref{lam1}).
To satisfy equations (\ref{lam2}), Zakharov proposed to impose
on the ``dressing matrix function'' $F_{ij} (s-u^i, s' -u^j)$ a certain
additional differential relation. If $F_{ij} (s - u^i, s' - u^j)$
satisfy the Zakharov differential relation, then
the rotation coefficients $\beta_{ij} (u)$ satisfy additionally
equations
(\ref{lam2}).

\bl
If both the function $F_{ij} (s-u^i, s' -u^j)$ and
the function
\be
\widetilde F_{ij} (s-u^i, s' -u^j) =
{\sqrt {f^j (u^j -s')} \over \sqrt {f^i (u^i -s)} }
F_{ij} (s-u^i, s' -u^j)   \label{f}
\ee
satisfy the Zakharov differential relation, then
the corresponding rotation coefficients $\beta_{ij} (u)$
(\ref{resh})
satisfy both equations (\ref{lam2}) and (\ref{lam3}).
\el

Actually, if $K_{ij} (s, s', u)$ is the solution of
the integral equation (\ref{int}) corresponding to the
function $F_{ij} (s-u^i, s' -u^j)$, then
\be
\widetilde K_{ij} (s, s', u) =
{\sqrt {f^j (u^j -s')} \over \sqrt {f^i (u^i -s)} }
K_{ij} (s, s', u)
\ee
is the solution of (\ref{int}) corresponding to
function (\ref{f}).
It is simple to prove multiplying the integral equation
(\ref{int}) by
$\sqrt {f^j (u^j -s')}/\sqrt {f^i (u^i -s)}: $

\be
\widetilde K_{ij}(s, s', u)  = \widetilde F_{ij}(s -u^i, s' - u^j)
+ \int_s^{\infty}
\sum_l \widetilde K_{il} (s, q, u)
\widetilde F_{lj} (q - u^l,  s'- u^j) dq.  \label{int2}
\ee

Then both $\widetilde \beta _{ij} (s, u) = \widetilde K_{ji} (s, s, u)$
and $\beta_{ij} (s, u) = K_{ji} (s, s, u)$
satisfy the Lam\'{e} equations (\ref{lam1}) and (\ref{lam2}).
Besides, we
have
\be
\widetilde \beta _{ij} (s, u) = \widetilde K_{ji} (s, s, u)
=
{\sqrt {f^i (u^i -s)} \over \sqrt {f^j (u^j -s)} }
K_{ji} (s, s, u) =
{\sqrt {f^i (u^i -s)} \over \sqrt {f^j (u^j -s)} }
\beta_{ij} (s, u).
\ee

Thus, in this case the rotation coefficients $\beta_{ij} (u)$
satisfy exactly all the equations (\ref{lam1})--(\ref{lam3}), that is,
they generate the corresponding compatible flat metrics.

The Zakharov differential reduction can be written as follows \cite{za}:

\be
{\pa F_{ij} (s, s', u) \over \pa s'} +
{\pa F_{ji} (s', s, u) \over \pa s} = 0. \label{za1}
\ee

Thus,
to resolve them for the matrix
function $F_{ij} (s- u^i, s' -u^j)$,
we can introduce $n(n-1)/2$ arbitrary functions of two variables
$\Phi (x, y),$ $i < j$, and put for $i<j$
\bea
&&
F_{ij} (s - u^i, s' - u^j) =
{\pa \Phi_{ij} (s-u^i, s' - u^j) \over \pa s},\nn\\
&&
F_{ji} (s - u^i, s' - u^j) =
- {\pa \Phi_{ij} (s' -u^i, s - u^j) \over \pa s},
\eea
and
\be
F_{ii} (s - u^i, s' - u^j) = {\pa \Phi_{ii} (s-u^i, s' - u^i) \over \pa s},
\ee
where $\Phi_{ii} (x, y)$, $i=1,...,N,$
are arbitrary skew-symmetric functions:
\be
\Phi_{ii} (x, y) = - \Phi_{ii} (y, x),
\ee
see \cite{za}.

For the function
\be
\widetilde F_{ij} (s-u^i, s' -u^j) =
{\sqrt {f^j (u^j -s')} \over \sqrt {f^i (u^i -s)} }
F_{ij} (s-u^i, s' -u^j)
\ee
the Zakharov differential relation (\ref{za1})
gives exactly $n(n-1)/2$ linear partial
 differential equations of the second order
 for
$n(n-1)/2$  functions $\Phi_{ij} (s- u^i, s' - u^j)$
of two variables:
\be
{\pa \over \pa s'} \left (
{\sqrt {f^j (u^j -s')} \over \sqrt {f^i (u^i -s)} }
{\pa \Phi_{ij} (s-u^i, s' - u^j) \over \pa s}
\right ) -
{\pa \over \pa s} \left (
{\sqrt {f^i (u^i -s)} \over \sqrt {f^j (u^j -s')} }
{\pa \Phi_{ij} (s-u^i, s' - u^j) \over \pa s'} \right )
\ee
or, equivalently,
\bea
&&
2 {\pa^2 \Phi_{ij} (s - u^i, s' - u^j) \over
\pa u^i \pa u^j} \left (f^i (u^i -s) - f^j (u^j - s') \right ) + \nn\\
&&
{\partial \Phi_{ij} (s- u^i, s' -u^j) \over \pa u^j}
{d f^i (u^i-s) \over d u^i}  -
{\partial \Phi_{ij} (s- u^i, s' -u^j) \over \pa u^i}
{d f^j (u^j-s') \over d u^j} =0.\label{ss}
\eea

It is interesting that all these equations (\ref{ss})
for functions $\Phi_{ij} (s- u^i, s' - u^j)$  are the same and
coincide with the corresponding single equation (\ref{lequa})
for the two-component case.

Besides, for $n$  functions
$\Phi_{ii} (s-u^i, s' - u^i) $ we have also $n$ linear partial
differential
equations of the second order from the Zakharov differential relation
(\ref{za1}):
\be
{\pa \over \pa s'} \left (
{\sqrt {f^i (u^i -s')} \over \sqrt {f^i (u^i -s)} }
{\pa \Phi_{ii} (s-u^i, s' - u^i) \over \pa s}
\right ) +
{\pa \over \pa s} \left (
{\sqrt {f^i (u^i -s)} \over \sqrt {f^i (u^i -s')} }
{\pa \Phi_{ii} (s' -u^i, s - u^i) \over \pa s'} \right )
\ee
or, equivalently,
\bea
&&
2 {\pa^2 \Phi_{ii} (s - u^i, s' - u^i) \over
\pa s \pa s'} \left (f^i (u^i -s) - f^i (u^i - s') \right ) - \nn\\
&&
{\partial \Phi_{ii} (s- u^i, s' -u^i) \over \pa s}
{d f^i (u^i-s') \over d s'}  +
{\partial \Phi_{ii} (s- u^i, s' -u^i) \over \pa s'}
{d f^i (u^i -s) \over d s} =0.  \label{ss2}
\eea


\medskip

\begin{flushleft}
Centre for Nonlinear Studies,\\
L.D.Landau Institute for Theoretical Physics,\\
Russian Academy of Sciences,\\
ul. Kosygina, 2,\\
Moscow, 117940  Russia\\
e-mail: mokhov@genesis.mi.ras.ru, mokhov@landau.ac.ru\\

\smallskip

Department of Mathematics,\\
University of Paderborn,\\
Paderborn, Germany\\
e-mail: mokhov@uni-paderborn.de
\end{flushleft}

\end{document}